# Improved Ratio Type Estimator Using Two Auxiliary Variables under Second Order Approximation

Prayas Sharma and Rajesh Singh

Department of Statistics, Banaras Hindu University

Varanasi-221005, India

Prayassharma02@gmail.com, rsinghstat@gmail.com

## Abstract

In this paper, we have proposed a new Ratio Type Estimator using auxiliary information on two auxiliary variables based on Simple random sampling without replacement (SRSWOR). The proposed estimator is found to be more efficient than the estimators constructed by Olkin (1958), Singh (1965), Lu (2010) and Singh and Kumar (2012) in terms of second order mean square error.

## 1. Introduction

In sampling survey, the use of auxiliary information is always useful in considerable reduction of the MSE of a ratio type estimator. Therefore, many authors suggested estimators using some known population parameters of an auxiliary variable. Hartley-Ross (1954), Quenouille's (1956) and Olkin (1958) have considered the problem of estimating the mean of a survey variable when auxiliary variables are made available. Jhajj and Srivastava (1983), Singh et al.(1995),

Upadhyaya and Singh (1999), Singh and Tailor (2003), Kadilar and Cingi (2006), Khoshnevisan et al. (2007), Singh et al. (2007), Singh and Kumar (2011,), etc. suggested estimators in simple random sampling using auxiliary variable.

Moreover, when two auxiliary variables are present Singh (1965,1967) and Perri(2007) suggested some ratio -cum -product type estimators. Most of these authors discussed the properties of estimators along with their first order bias and MSE. Hossain et al. (2006), Sharma et al. (2013a, b) studied the properties of some estimators under second order approximation. In this paper, we have suggested an estimator using auxiliary information in simple random sampling and compared with some existing estimators under second order of approximation when information on two auxiliary variables are available.

Let U= ($U_1$ ,$U_2$ , $U_3$, …..,$U_i$, ….$U_N$ ) denotes a finite population of distinct and identifiable units. For estimating the population mean $\overline{Y}$ of a study variable Y, let us consider X and Z are the two auxiliary variable that are correlated with study variable Y, taking the corresponding values of the units. Let a sample of size n be drawn from this population using simple random sampling without replacement (SRSWOR) and $y_i$ , $x_i$ and $z_i$ (i=1,2,…..n ) are the values of the study variable and auxiliary variables respectively for the i-th units of the sample.

When the information on two auxiliary variables are available Singh (1965,1967) proposed some ratio-cum-product estimators in simple random sampling without replacement to estimate the population mean $\overline{Y}$ of the study variable y, generalized version of one of these estimators is given by,

$$t_1 = \bar{y}\left(\frac{\overline{X}}{\bar{x}}\right)^{\alpha_1}\left(\frac{\overline{Z}}{\bar{z}}\right)^{\alpha_2} \tag{1.1}$$

where $\alpha_1$ and $\alpha_2$ are suitably chosen scalars such that the mean square error of $t_1$ is minimum.

and $\bar{y} = \frac{1}{n}\sum_{i=1}^{n} y_i$ , $\bar{x} = \frac{1}{n}\sum_{i=1}^{n} x_i$ and $\bar{z} = \frac{1}{n}\sum_{i=1}^{n} z_i$ ,

Olkin(1958) proposed an estimator $t_2$ as-

$$t_2 = \bar{y}\left[\lambda_1 \frac{\bar{X}}{\bar{x}} + \lambda_2 \frac{\bar{Z}}{\bar{z}}\right] \tag{1.2}$$

where, $\lambda_1$ and $\lambda_2$ are the weights that satisfy the condition $\lambda_1 + \lambda_2 = 1$

Lu (2010) proposed multivariate ratio estimator using information on two auxiliary variables as-

$$t_3 = \bar{y}\left[\frac{w_1\bar{X}_1 + w_2\bar{X}_2}{w_1\bar{x}_1 + w_2\bar{x}_2}\right]^{\alpha} \tag{1.3}$$

where, $w_1$ and $w_2$ are weights that satisfy the condition: $w_1 + w_2 = 1$.

Singh and Kumar (2012) proposed exponential ratio-cum-product estimator. Generalized form of this estimator is given by

$$t_4 = \bar{y}\exp\left[\frac{\bar{X} - \bar{x}}{\bar{X} + \bar{x}}\right]^{\beta_1} \exp\left[\frac{\bar{Z} - \bar{z}}{\bar{Z} + \bar{z}}\right]^{\beta_2} \tag{1.4}$$

where $\beta_1$ and $\beta_2$ are constants such that the MSE of estimator $t_4$ is minimise.

**Theorem1.1** Let $e_0 = \frac{\bar{y} - \bar{Y}}{\bar{Y}}$ , $e_1 = \frac{\bar{x} - \bar{X}}{\bar{X}}$ and $e_2 = \frac{\bar{z} - \bar{Z}}{\bar{Z}}$

then $E(e_0) = E(e_1) = E(e_2) = 0$ and variances and co-variances are as follows:

(i) $V(e_0) = E\{(e_0)^2\} = \frac{L_1 C_{200}}{\bar{Y}^2} = V_{200}$

(ii) $V(e_1) = E\{(e_1)^2\} = \frac{L_1 C_{020}}{\bar{X}^2} = V_{020}$

(iii) $V(e_1) = E\{(e_2)^2\} = \frac{L_1 C_{002}}{\bar{Z}^2} = V_{002}$

(iv) $COV(e_0, e_1) = E\{(e_0 e_1)\} = \frac{L_1 C_{110}}{\bar{X}\bar{Y}} = V_{110}$

(v) $COV(e_1, e_2) = E\{(e_1 e_2)\} = \frac{L_1 C_{011}}{\bar{X}\bar{Z}} = V_{011}$

(vi) $COV(e_0, e_2) = E\{(e_0 e_2)\} = \frac{L_1 C_{101}}{\bar{Y}\bar{Z}} = V_{101}$

(vii) $E\{(e_0{}^2 e_1)\} = \frac{L_2 C_{210}}{\bar{Y}^2 \bar{X}} = V_{210}$

(viii) $E\{(e_0{}^2 e_2)\} = \frac{L_2 C_{201}}{\bar{Y}^2 \bar{Z}} = V_{201}$

(ix) $E\{(e_1{}^2 e_2)\} = \frac{L_2 C_{021}}{\bar{X}^2 \bar{Z}} = V_{021}$

(x) $E\{(e_0 e_1^2)\} = \frac{L_2 C_{120}}{\bar{Y}\bar{X}^2} = V_{120}$

(xi) $E\{(e_1 e_2^2)\} = \frac{L_2 C_{012}}{\bar{X}\bar{Z}^2} = V_{012}$

(xii) $E\{(e_0 e_2^2)\} = \frac{L_2 C_{102}}{\bar{Y}\bar{Z}^2} = V_{102}$

(xiii) $E\{(e_1{}^3)\} = \frac{L_2 C_{030}}{\bar{X}^3} = V_{030}$

(xiv) $E(e_1{}^3 e_2) = \frac{L_3 C_{031} + 3L_4 C_{020} C_{011}}{\bar{X}^3 \bar{Z}} = V_{031}$

(xv) $E(e_1 e_2^3) = \frac{L_3 C_{013} + 3L_4 C_{002} C_{011}}{\bar{X}\bar{Z}^3} = V_{013}$

(xvi) $E(e_0 e_1{}^3) = \frac{L_3 C_{130} + 3L_4 C_{020} C_{110}}{\bar{Y}\bar{X}^3} = V_{130}$

where, $L_1 = \frac{(N-n)}{(N-1)}\frac{1}{n}$, $L_2 = \frac{(N-n))(N-2n)}{(N-1)(N-2)}\frac{1}{n^2}$

$$L_3 = \frac{(N-n))(N^2+N-6nN+6n^2)}{(N-1)(N-2)(N-3)}\frac{1}{n^3}, \quad L_4 = \frac{N(N-n))(N-n-1)(n-1)}{(N-1)(N-2)(N-3)}\frac{1}{n^3}$$

And $$Cpqr = \sum_{i=1}^{N} (Xi - \overline{X})^p (Y_i - \overline{Y})^q \left(Z_i - \overline{Z}\right)^r$$

This theorem will be used to obtain MSE expressions of estimators considered here. Proof of this theorem is straight forward by using SRSWOR ( see Sukhatme and Sukhatme (1970)).

## 2. First Order Biases and Mean Squared Errors

The expression of the biases of the estimators $t_{1,}$ $t_2$ $t_3$ and $t_4$ to the first order of approximation are respectively, written as

$$Bias(t_1) = \overline{Y}\left[-\alpha_1 V_{110} + \alpha_2 V_{102} + R_1 V_{020} + S_1 V_{002} + \alpha_1\alpha_2 V_{011}\right] \tag{2.1}$$

where , $R_1 = \frac{\alpha_1(\alpha_1+1)}{2}$ $S_1 = \frac{\alpha_2(\alpha_2+1)}{2}$.

$$Bias(t_2) = \overline{Y}\left[-w_1 V_{110} - w_2 V_{101} + w_1 V_{020} + w_2 V_{002} + \alpha_1\alpha_2 V_{011}\right] \tag{2.2}$$

$$Bias(t_3) = \overline{Y}\Big[-\alpha\lambda\left(w_1\overline{X}_1 V_{110} - w_2 V_{101}\right) + \lambda^2\alpha\frac{(\alpha+1)}{2}\left(w_1^2\overline{X}_1^2 V_{020} + w_2^2\overline{X}_2^2 V_{002} + 2w_1 w_1 \overline{X}_1\overline{X}_2 V_{011}\right)\Big] \tag{2.3}$$

where, $\lambda = \frac{1}{w_1\overline{X}_1 + w_2\overline{X}_2}$

$$Bias(t_4) = \overline{Y}\left[-\frac{\beta_1}{2}V_{110} - \frac{\beta_2}{2}V_{101} + V_{020}\left(\frac{\beta_1}{4} + \frac{\beta_1^2}{8}\right) + \left(\frac{\beta_2}{4} + \frac{\beta_2^2}{8}\right)V_{002}\right] \tag{2.4}$$

Expressions for the MSE of the estimators $t_1$, $t_2$, $t_3$ and $t_4$ to the first order of approximation are respectively, given by

$$MSE(t_1) = \bar{Y}^2 \left[ V_{200} + \alpha_1^{\,2} V_{020} + \alpha_2^2 V_{002} - 2\alpha_1 V_{110} - 2\alpha_2 V_{101} + 2\alpha_1\alpha_2 V_{011} \right] \quad (2.5)$$

The MSE of the estimator $t_1$ is minimized for

$$\alpha_1^* = \left[ \frac{\rho_{yx} - \rho_{yz}\rho_{xz}}{1-\rho_{xz}^2} \right] \sqrt{\frac{V_{200}}{V_{020}}} \quad (2.6)$$

And

$$\alpha_1^* = \left[ \frac{\rho_{yz} - \rho_{yx}\rho_{xz}}{1-\rho_{xz}^2} \right] \sqrt{\frac{V_{200}}{V_{002}}} \quad (2.7)$$

where $\alpha_1^*$ and $\alpha_2^*$ are, respectively, partial regression coefficients of y on x and of y on z in simple random sampling.

$$MSE(t_2) = \bar{Y}^2 \left[ V_{200} + \lambda_1^{\,2} V_{020} + \lambda_2^2 V_{002} - 2\lambda_1 V_{110} - 2\lambda_2 V_{101} + 2\lambda_1\lambda_2 V_{011} \right] \quad (2.8)$$

The MSE of the estimator $t_2$ is minimum for

$$\lambda_1^* = \frac{V_{002} - V_{101} + V_{110} - V_{012}}{V_{020} + V_{002} - 2V_{012}} \quad \text{and } \lambda_2^* = 1 - \lambda_1^* \quad (2.9)$$

$$MSE(t_3) = \bar{Y}^2 \left[ V_{200} + \alpha^2\lambda^2 (w_1^2 \bar{X}_1^2 V_{020} + w_2^2 \bar{X}_2^2 V_{002} + 2 w_1 w_2 \bar{X}_1 \bar{X}_2 V_{011}) \right.$$

$$\left. - 2\alpha\lambda (w_1 \bar{X}_1 V_{110} + w_2 \bar{X}_2 V_{101} \right] \quad (2.10)$$

Differentiating (2.10) with respect to $w_1$ and $w_2$ partially, we get the optimum values of $w_1$ and $w_2$ respectively as

$$w_1^* = \frac{\bar{X}_1 V_{110} - \bar{X}_2 V_{101} + \bar{X}_2^2 V_{002} - \bar{X}_1 \bar{X}_2 V_{011}}{\bar{X}_1^2 V_{020} + \alpha\lambda \bar{X}_2^2 V_{002} - 2\bar{X}_1 \bar{X}_2 V_{011}} \quad \text{and } w_2^* = 1 - w_1^* \quad (2.11)$$

For optimum value of $w_1 = w_1^*$ and $w_2 = w_2^*$, MSE of the estimator $t_3$ is minimum.

$$MSE(t_4) = \bar{Y}^2\left[V_{200} + \frac{\beta_1^{\ 2}}{4}V_{020} + \frac{\beta_2^{\ 2}}{4}V_{002} - \beta_1 V_{110} - \beta_2 V_{101} + \frac{\beta_1\beta_2}{2}V_{011}\right] \tag{2.12}$$

On differentiating (2.12) with respect to $\beta_1$ and $\beta_2$ respectively, we get the optimum values of $\beta_1$ and $\beta_2$ as

$$\beta_1^* = \frac{2(V_{110}V_{002} - V_{101}V_{011})}{(V_{002}V_{020} - V_{011}^2)} \tag{2.13}$$

and

$$\beta_2^* = \frac{2(V_{020}V_{101} - V_{110}V_{011})}{(V_{002}V_{020} - V_{011}^2)} \tag{2.14}$$

Estimators $t_1$, $t_2$ $t_3$ and $t_4$ at their respective optimum values attains MSE values which are equal to the MSE of regression estimator for two auxiliary variables.

## 3. Proposed Estimator

When auxiliary information on two auxiliary variables are known, we propose an estimator $t_5$ as

$$t_5 = \bar{y}\left[k_1\left\{\frac{c\bar{X} - d\bar{x}}{(c-d)\bar{X}}\right\}^{\delta_1} + k_2\left\{2 - \left(\frac{\bar{z}}{\bar{Z}}\right)^{\delta_2}\right\}\right] \tag{3.1}$$

where d and c are either real numbers or a function of the known parameters associated with auxiliary information . $k_1$ and $k_2$ are constants to be determined such that the MSE of estimator $t_5$ is minimum under the condition that $k_1+k_2 = 1$ and $\delta_1$ and $\delta_2$ are integers and can take values -1, 0 and +1 .

Expressing the estimator $t_5$ in terms of e's we have

$$t_5 = \bar{Y}\left[k_1(1-\eta_1 e_1)^{\delta_1} + k_2\left\{2-(1+e_2)^{\delta_2}\right\}\right] \tag{3.2}$$

where, $\eta_1 = \dfrac{d}{c-d}$.

We assume that $|\eta_1 e_1| < 1$ so that $(1-\eta_1 e_1)^{\delta_1}$ are expandable. Expanding the right hand side of (3.2), and neglecting terms of *e*'s having power greater than two we have

$$\left(t_5 - \bar{Y}\right) = \bar{Y}\left[-k_1\eta_1\delta_1 e_0 e_1 - k_2\delta_2 e_0 e_2 + k_1 M_1 e_1^2 - k_2 N_1 e_2^2 + \alpha_1\alpha_2 e_1 e_2\right] \tag{3.3}$$

Taking expectation on both sides we get bias of estimator $t_5$, to the first degree of approximation as

$$\text{Bias}(t_5) = \bar{Y}\left[-k_1\eta_1\delta_1 V_{110} - k_2\delta_2 V_{101} + k_1 M_1 V_{020} - k_2 N_1 V_{002} + \alpha_1\alpha_2 V_{011}\right] \tag{3.4}$$

where,

$$M_1 = \frac{\delta_1(\delta_1 - 1)}{2}\eta_1^2, \qquad N_1 = \frac{\delta_2(\delta_2 - 1)}{2}.$$

Squaring both sides of (3.3) and neglecting terms of *e*'s having power greater than two we have

$$(t_5 - \bar{Y})^2 = \bar{Y}\left[e_0^2 + k_1^2\eta_1^2\delta_1^2 e_1^2 + k_2^2\delta_2^2 e_2^2 - 2k_1\eta_1\delta_1 e_0 e_1 - 2k_2\delta_2 e_0 e_2 + 2k_1 k_2\delta_1\delta_2\eta_1 e_1 e_2\right] \tag{3.5}$$

Taking expectation on both sides of (3.5) and using theorem 1.1, we get the MSE of $t_5$ up to first degree of approximation as

$$\text{MSE}(t_5) = \bar{Y}^2\left[V_{200} + k_1^2\delta_1^2\eta_1{}^2 V_{020} + k_2^2\delta_2^2 V_{002} - 2k_1\delta_1\eta_1 V_{110} - 2k_2\delta_2 V_{101} + 2k_1 k_2\delta_1\delta_2\eta_1 V_{011}\right] \tag{3.6}$$

Differentiating (3.6) with respect to $k_1$ and $k_2$ partially, equating them to zero and after simplification, we get the optimum values of $k_1$ and $k_2$ respectively, as

$$k_1^* = \frac{\delta_1\eta_1 V_{110} + \delta_2^2 V_{002} - \delta_2 V_{102} - \delta_1\delta_2\eta_1 V_{011}}{\delta_1^2\eta_1^2 V_{110} + \delta_2^2 V_{002} - 2\delta_1\delta_2\eta_1 V_{011}}, \qquad k_2^* = 1 - k_1^*. \tag{3.7}$$

Putting these values in (3.6) we get minimum MSE of estimator $t_5$. The minimum MSE of the estimator $t_1$, $t_2$, $t_3$, $t_4$ and proposed estimator $t_5$ is equal to the MSE of combined regression estimator based on two auxiliary variables, which motivated us to study the properties of estimators up to the second order of approximation.

## 4. Second Order Biases and Mean Squared Errors

Expressing estimator $t_i$'s (i=1,2,3,4,5) in terms of e's (i=0,1), we get

$$t_1 = \overline{Y}(1+e_0)\left\{(1+e_1)^{-\alpha_1}(1+e_1)^{-\alpha_2}\right\} \tag{4.1}$$

Or

$$\begin{aligned}(t_1-\overline{Y}) = \overline{Y}\{&e_0 - \alpha_1 e_1 - \alpha_1 e_0 e_1 - \alpha_2 e_2 - \alpha_2 e_0 e_2 + \alpha_1\alpha_2 e_1 e_2 + \alpha_1\alpha_2 e_0 e_1 e_2 + R_1 e_1{}^2 + R_1 e_0 e_1{}^2 \\ &+ S_1 e_2^2 + S_1 e_0 e_2^2 - R_2 e_1^3 - R_2 e_0 e_1^3 - S_2 e_2^3 - S_2 e_0 e_2^3 - R_2 e_1^3 - \alpha_2 R_1 e_2 e_1^2 - \alpha_1 S_1 e_1 e_2^2 \\ &+ \alpha_2 R_2 e_1^3 e_2 + \alpha_1 R_1 e_1 e_2^3\}\end{aligned} \tag{4.2}$$

Squaring both sides and neglecting terms of *e*'s having power greater than four, we have

$$\begin{aligned}(t_1-\overline{Y})^2 = \overline{Y}^2\big[&e_0^2 + \alpha_1^2 e_1^2 + \alpha_2^2 e_2^2 - 2\alpha_1 e_0 e_1 - 2\alpha_2 e_0 e_2 + 2\alpha_1\alpha_2 e_1 e_2 - 2\alpha_1 e_0^2 e_1 - 2\alpha_2 e_0^2 e_2 \\ &+ (2R_1 + 2\alpha_1^2) e_0 e_1^2 + 2S_1 e_0 e_2^2 - 2\alpha_1^2\alpha_2 e_1^2 e_2 - 2S_1\alpha_1 e_1 e_2^2 - 2R_1\alpha_1 e_1^3 + 6\alpha_1\alpha_2 e_0 e_1 e_2 \\ &+ (\alpha_1^2 + 2R_1) e_0^2 e_1^2 + (\alpha_2^2 + 2S_1) e_0^2 e_2^2 + (\alpha_1^2\alpha_2^2 + 2R_1 S_1) e_1^2 e_2^2 - (4\alpha_1^2\alpha_2 + 6R_1\alpha_1) e_0 e_1^2 e_2 \\ &- 4\alpha_1(S_1 + \alpha_2^2) e_0 e_1 e_2^2 + 4\alpha_1\alpha_2 e_0^2 e_1 e_2 - 2(R_2 + 2\alpha_1 R_1) e_0 e_1^3 - 2(S_2 + 2\alpha_2 S_1) e_0 e_2^3 \\ &- 2\alpha_1(S_2 - \alpha_2 S_1) e_1 e_2^2 + 2\alpha_2(R_2 + \alpha_1 R_1) e_1^2 e_2 + (R_1^2 + 2\alpha_1 R_2) e_1^4 + (S_1^2 + 2\alpha_2 S_2) e_2^4\big]\end{aligned} \tag{4.3}$$

Taking expectations, and using theorem 1.1 we get the MSE of the estimator $t_1$ up to the second order of approximation as

$$\begin{aligned}MSE_2(t_1) = \overline{Y}^2\big[&V_{200} + \alpha_1^2 V_{020} + \alpha_2^2 V_{002} - 2\alpha_1 V_{110} - 2\alpha_2 V_{101} + 2\alpha_1\alpha_2 V_{011} - 2\alpha_1 V_{210} - 2\alpha_2 V_{201} \\ &+ (2R_1 + 2\alpha_1^2) V_{120} + 2S_1 V_{102} - 2\alpha_1^2\alpha_2 V_{021} - 2S_1\alpha_1 V_{012} - 2R_1\alpha_1 V_{030} + 6\alpha_1\alpha_2 V_{111}\end{aligned}$$

$$+(\alpha_1^2+2R_1)V_{220}+(\alpha_2^2+2S_1)V_{202}+(\alpha_1^2\alpha_2^2+2R_1S_1)V_{022}-(4\alpha_1^2\alpha_2+6M_1\alpha_1)V_{121}$$

$$-4\alpha_1(S_1+\alpha_2^2)V_{112}+4\alpha_1\alpha_2V_{211}-2(R_2+2\alpha_1R_1)V_{130}-2(S_2+2\alpha_2S_1)V_{103}$$

$$-2\alpha_1(S_2-\alpha_2S_1)V_{012}+2\alpha_2(R_2+\alpha_1R_1)V_{021}+(R_1^2+2\alpha_1R_2)V_{040}+(S_1^2+2\alpha_2S_2)V_{004}\Big] \quad (4.4)$$

where, $R_1=\frac{\alpha_1(\alpha_1+1)}{2}$, $R_2=\frac{\alpha_1(\alpha_1+1)(\alpha_1+2)}{6}$, $R_3=\frac{\alpha_1(\alpha_1+1)(\alpha_1+2)(\alpha_1+3)}{24}$,

$S_1=\frac{\alpha_2(\alpha_2+1)}{2}$, $S_2=\frac{\alpha_2(\alpha_2+1)(\alpha_2+2)}{6}$, $S_3=\frac{\alpha_2(\alpha_2+1)(\alpha_2+2)(\alpha_2+3)}{24}$.

Similarly, MSE expression of estimator $t_2$ is given by

$$MSE_2(t_2)=\overline{Y}^2\Big[V_{200}+\lambda_1^2\left(V_{020}+V_{220}+3V_{040}+2V_{120}-2V_{030}-4V_{130}\right)$$

$$+\lambda_2^2\left(V_{002}+V_{202}+3V_{004}+2V_{102}-2V_{003}-4V_{103}\right)$$

$$+2\lambda_1\left(-V_{110}-V_{210}+V_{120}+V_{220}-V_{130}\right)+2\lambda_2\left(-V_{101}-V_{201}+V_{202}\right)$$

$$+2\lambda_1\lambda_2\left(V_{011}+2V_{111}-V_{012}-2V_{112}+V_{013}+V_{211}-V_{021}-2V_{121}+V_{022}+V_{031}\right)\Big] \quad (4.5)$$

MSE expression of estimator $t_3$ is given by

$$MSE_2(t_3)=\overline{Y}^2\Big[V_{200}+w_1^2\overline{X}_1\left\{\alpha^2\theta^2(V_{020}+V_{220}+2V_{120})+2A_1\theta^2(V_{120}+V_{220})\right\}$$

$$+w_2^2\overline{X}_2\left\{\alpha^2\theta^2(V_{002}+V_{202}+2V_{120})+2A_1\theta^2(V_{120}+V_{202})\right\}$$

$$+2w_1w_2\overline{X}_1\overline{X}_2\left\{\alpha^2\theta^2(V_{211}+2V_{111})+2A_1\theta^2(V_{111}+V_{211}\right\}$$

$$+2w_1^3w_2\overline{X}_1^3\overline{X}_2\left(A_1^2\theta^4+4A_2\alpha\theta^4\right)V_{031}+2w_1w_2^3\overline{X}_1\overline{X}_2^3\left(A_1^2\theta^4+4A_2\alpha\theta^4\right)V_{013}$$

$$-2\alpha\theta w_1\overline{X}_1\left(V_{110}+V_{210}\right)-2\alpha\theta w_2\overline{X}_2\left(V_{101}+V_{201}\right)$$

$$+3w_1^2w_2\overline{X}_1^2\overline{X}_2\left(-2A_2\theta^3V_{121}-4A_1\alpha\theta^3V_{121}-2A_1\alpha\theta^3V_{012}\right)$$

$$+3w_1\overline{X}_1w_2^2\overline{X}_2^2\left(-2A_2\theta^3V_{112}-4A_1\alpha\theta^3V_{112}-2A_1\alpha\theta^3V_{012}\right)$$

$$+w_1^4\overline{X}_1^4\left(A_1^4\theta^4+2A_2\alpha\theta^4\right)V_{040}+w_2^4\overline{X}_2^4\left(A_1^4\theta^4+2A_2\alpha\theta^4\right)V_{002}$$

$$+w_1^3\overline{X}_1^3\left(-2A_2\theta^3V_{130}-4A_1\alpha\theta^3V_{130}-2A_1\alpha\theta^3V_{030}\right)$$

$$+w_2^3\overline{X}_2^3\left(-2A_2\theta^3V_{103}-4A_1\alpha\theta^3V_{103}-2A_1\alpha\theta^3V_{003}\right)$$

$$+6w_1^2\overline{X}_1^2w_2^2\overline{X}_2^2\left(A_1^2\theta^4+2A_2\alpha\theta^4\right)V_{022}\Big] \tag{4.6}$$

MSE expression of estimator $t_4$ is given by

$$MSE_2(t_4)=\overline{Y}^2\Big[V_{200}+\frac{\beta_1^2}{4}V_{020}+\frac{\beta_2^2}{4}V_{002}-\beta_1V_{110}-\beta_1V_{210}\left\{M+\frac{\beta_1^2}{4}\right\}+\frac{V_{102}}{2}\left\{N+\frac{\beta_2^2}{4}\right\}$$

$$-V_{021}\left\{M+\beta_1^2\beta_2\right\}-V_{012}\left\{N+\beta_2^2\beta_1\right\}+\frac{3}{2}\beta_1\beta_2V_{111}+\frac{1}{2}\beta_1\beta_2V_{011}-\frac{3}{4}\beta_1MV_{030}-\frac{1}{2}\beta_2NV_{003}$$

$$+\frac{V_{220}}{2}\left\{M+\frac{\beta_1^2}{4}\right\}+\frac{V_{202}}{2}\left\{N+\frac{\beta_2^2}{4}\right\}+\frac{V_{022}}{8}\left\{\frac{\beta_1^2\beta_2^2}{2}+\beta_1+\beta_2Q+MN\right\}-\frac{V_{130}}{4}\left\{O+2\beta_1M\right\}$$

$$-\frac{V_{103}}{4}\left\{P+2\beta_2N\right\}+\frac{V_{031}}{8}\left\{\beta_1Q+O+SM\right\}+\frac{V_{013}}{8}\left\{\beta_2R+\beta_1P+SN\right\}$$

$$-\frac{V_{121}}{4}\left\{Q+2\beta_1^2\beta_2+\beta_2M\right\}-\frac{V_{112}}{4}\left\{2\beta_2^2\beta_1+2\beta_1N+R\right\}+-\frac{V_{040}}{16}\left\{2\beta_1O+M^2\right\}$$

$$+\frac{V_{004}}{16}\left\{2\beta_2P+N^2\right\}+\frac{V_{211}}{2}\left\{\beta_1^2\beta_2+S\right\}\Bigg] \tag{4.7}$$

where,

$$M=\left(\beta_1+\frac{\beta_1^2}{2}\right),\quad N=\left(\beta_2+\frac{\beta_2^2}{2}\right),\quad O=\left(\beta_1^2+\frac{\beta_1^3}{6}\right)$$

$$P=\left(\beta_2^2+\frac{\beta_2^3}{6}\right),\quad Q=\left(\beta_1\beta_2+\frac{\beta_1^2\beta_2}{2}\right),\quad R=\left(\beta_1\beta_2+\frac{\beta_1\beta_2^2}{2}\right)$$

**MSE of Proposed estimator $t_5$ up to second order of approximation**

Expanding right hand side of equation (3.2) and neglecting terms of *e*'s having power greater than four, we get

$$t_5 - \overline{Y} = \overline{Y}\Big[e_0 + k_1\left(-\delta_1\eta_1 e_1 + M_1 e_1^2 - M_2 e_1^3 + M_3 e_1^4\right) + k_2\left(-\delta_2 e_2 + N_1 e_2^2 - N_2 e_2^3 + N_3 e_2^4\right)$$

$$+ k_1\left(-\delta_1\eta_1 e_0 e_1 + M_1 e_0 e_1^2 - M_2 e_0 e_1^3\right) + k_2\left(-\delta_2 e_0 e_2 + N_1 e_0 e_2^2 - N_2 e_0 e_2^3\right)\Big] \quad (4.8)$$

Squaring both sides of (4.8) and neglecting terms of *e*'s having power greater than four, we have

$$\left(t_5 - \overline{Y}\right)^2 = \overline{Y}^2\Big[e_0^2 + k_1^2\left\{\delta_1^2\eta_1^2(e_1^2 + e_1^2 e_2^2 + 2e_0 e_1^2) + 2\delta_1\eta_1(M_1 e_1^4 - 2M_1 e_1^3) + M_1^2 e_1^4\right\}$$

$$+ k_2^2\left\{\delta_2^2(e_2^2 + e_0^2 e_2^2 + 2e_0 e_2^2) + 2\delta_2(N_1 e_0 e_2^3 + N_1 e_1^2 + N_2 e_2^4) + N_1^2 e_2^4\right\}$$

$$+ 2k_1\left\{\delta_1\eta_1(-e_0 e_1 - e_0^2 e_1) + M_1(e_0 e_1^2 + e_0^2 e_1^2) - M_2 e_0 e_1^3\right\}$$

$$+ 2k_2\left\{\delta_2(-e_0 e_2 - e_0^2 e_2) + N_1(-e_0 e_2^2 - e_0^2 e_2^2) - N_2 e_0 e_2^3\right\}$$

$$+ 2k_1 k_2\left\{\delta_1\delta_2\eta_1(e_1 e_2 + 2e_0 e_1 e_2 + e_0^2 e_1 e_2) + \delta_1\eta_1(N_1 e_1 e_2^2 + N_2 e_1 e_2^3)\right\}\Big] \quad (4.9)$$

Taking expectations on both sides of (4.9) and using theorem1.1, we get the MSE of estimator $t_5$,up to the second order of approximation as

$$MSE_2(t_5) = \overline{Y}^2\Big[V_{200} + k_1^2\left\{\delta_1^2\eta_1^2(V_{020} + V_{022} + 2V_{120}) + 2\delta_1\eta_1(M_1 V_{040} - 2M_1 V_{030}) + M_1^2 V_{040}\right\}$$

$$+ k_2^2\left\{\delta_2^2(V_{002} + V_{202} + 2V_{102}) + 2\delta_2(N_1 V_{103} + N_1 V_{030} + N_2 V_{004}) + N_1^2 V_{004}\right\}$$

$$+ 2k_1\left\{\delta_1\eta_1(-V_{110} - V_{210}) + M_1(V_{120} + V_{220}) - M_2 V_{130}\right\}$$

$$+ 2k_2\left\{\delta_2(-V_{101} - V_{201}) + N_1(-V_{102} - V_{202}) - N_2 V_{103}\right\}$$

$$+ 2k_1 k_2\left\{\delta_1\delta_2\eta_1(V_{011} + 2V_{111} + V_{211}) + \delta_1\eta_1(N_1 V_{012} + N_2 V_{013})\right\}\Big] \quad (4.10)$$

## 5. Numerical Illustration

For a natural population data, we have calculated the mean square error's of the estimator's and compared MSE's of the estimator's under first and second order of approximations.

**Data Set**

The data for the empirical analysis are taken from Book, "An Introduction to Multivariate Statistical Analysis", page no. 58, 2nd Edition By T.W. Anderson.

The population consist of 25 persons with Y= Head length of second son, X= Head length of first son and Z= Head breadth of first son. The following values are obtained from Raw data given on page no. 58.

$\overline{Y} = 183.84,\ \overline{X} = 185.72,\ \overline{Z} = 151.12, N = 25, n = 7$ and

V020=0.000244833,V200=0.000306792,V020=0.000284171,V110=0.00020986,V011=0.000193753,
V101= 0.000189972, V111= -0.000059732, V210 = -0.0000004582, V201= -0.0000002.77,
V021= -0.0000002775, V102= -0.0000002354, V120= -0.00000036455, V012=0.00000025179,
V202=0.000013,V003=0.000001544,V031=0.3893411 ,V013=0.380025,
V030= -0.00001363,V103=0.00001215,V040=0.0000214,V220=0.000015
V022=0.001624,V121=0.0000115,V211=0.0000122,V004=0.00001384,V112=0.000000085

**Table 5.1: MSE's of the estimators $t_i$ (i=1,2,3,4,5)**

| Estimators | MSE | |
|---|---|---|
| | First order | Second order |
| $t_1$ | 4.508 | 16156.644 |
| $t_2$ | 4.508 | 27204.321 |
| $t_3$ | 4.508 | 17679.890 |
| $t_4$ | 4.508 | 20928.689 |
| $t_5$ | 4.508 | 275.926* |

*for $\delta_1 = \delta_2 = 1$

This table shows the comparison of estimators on the basis of MSE because $\bar{y}$ cannot be extended up to second order of approximation therefore, we are unable to calculate PRE for second order approximation.

## 6. Conclusion

In the Table 5.1 the MSE's of the estimators $t_1$, $t_2$ , $t_3$, $t_4$ and $t_5$ are written under first order and second order of approximations. It has been observed that for all the estimators, the mean squared error increases for second order. Observing second order MSE's we conclude that the estimator $t_5$ is best estimator among the estimators considered here for the given data set.

Acknowledgement: The authors are thankful to the learned referee's for their valuable suggestions leading to the improvement of contents and presentation of the original manuscript.

**REFERENCES**

1. Anderson, T. W. (1984). An introduction to multivariate statistical analysis (2nd ed.). New York: Wiley.
2. Hartley, H.O. and Ross,A. (1954): Unbiased ratio estimators. **Nature**, 174, 270-271.
3. Hossain, M.I., Rahman, M.I. and Tareq, M.(2006) : Second order biases and mean squared errors of some estimators using auxiliary variable. SSRN.
4. Jhajj H.S and Srivastva S.K., (1983): A class of PPS estimators of population mean using auxiliary information. *Jour. Indian Soc. Agril. Statist.* 35, 57-61.
5. Kadilar, C.; Cingi, H., (2006): Improvement in Estimating the Population Mean in Simple Random Sampling, Applied Mathematics Letters, 19, 1, 75-79.
6. Khoshnevisan, M., Singh, R., Chauhan, P., Sawan, N., and Smarandache, F. (2007). A general family of estimators for estimating population mean using known value of some population parameter(s), Far East Journal of Theoretical Statistics 22 181–191.

7. Lu,J. Yan,Z. Ding,C. Hong,Z., (2010): 2010 International Conference On Computer and Communication Technologies in Agriculture Engineering (CCTAE), 3, 136-139.

8. Olkin, I. (1958) Multivariate ratio estimation for finite populations, *Biometrika*, 45, 154–165.

9. Perri, P. F. (2007). Improved ratio-cum-product type estimators. *Statist. Trans.* 8:51–69.

10. Quenouille, M. H. (1956): Notes on bias in estimation, *Biometrika*, 43, 353-360.

11. Sharma, P.**,** Singh, R. ,Jong, Min-kim. (2013a): Study of Some Improved Ratio Type Estimators using information on auxiliary attributes under second order approximation. Journal of Scientific Research, vol.57, 138-146.

12. Sharma, P., Singh, R.,Verma, H. **(**2013b): Some Exponential Ratio-Product Type Estimators using information on Auxiliary Attributes under Second Order Approximation. International Journal of Statistics and Economics, Vol. 12, issue no. 3.

13. Singh, H.P. and Tailor, R. (2003). Use of known correlation coefficient in estimating the finite population mean. Statistics in Transition 6 555-560.

14. Singh, M. P. (1965). On the estimation of ratio and product of the population parameters. Sankhya *B* 27:231-328.

15. Singh, M. P. (1967). Ratio cum product method of estimation. Metrika 12:34–42.

16. Singh, R. and Kumar, M. (2011): A note on transformations on auxiliary variable in survey sampling. MASA, 6:1, 17-19.

17. Singh, R. and Kumar, M. (2012 Improved Estimators of Population Mean Using Two Auxiliary Variables in Stratified Random Sampling Pak.j.stat.oper.res. Vol.VIII No.1 pp65-72.

18. Singh, R., Cauhan, P., Sawan, N., and Smarandache, F. (2007). Auxiliary Information and A Priori Values in Construction of Improved Estimators. Renaissance High Press.

19. Singh S, Mangat N.S and Mahajan P.K., (1995): General Class of Estimators. *Jour. Indian Soc. Agril. Statist.*47(2), 129-133.

20. Srivastava, S.K. (1967) : An estimator using auxiliary information in sample surveys. Cal. Stat. Ass. Bull. 15:127-134.

21. Sukhatme, P.V. and Sukhatme, B.V. (1970): Sampling theory of surveys with applications. Iowa State University Press, Ames, U.S.A.

22. Upadhyaya, L. N. and Singh, H. P. (1999): Use of transformed auxiliary variable in estimating the finite population mean. Biom. Jour., 41, 627-636.